\documentclass[twoside,11pt,leqno]{article}
\usepackage{amsfonts}

 \def\X{{\cal X}}  \def\H{{\cal H}}

\def\B{B({\cal H})} \def\b{B({\cal X})}

\newtheorem{df}{Definition}[section]
\newtheorem{thm}[df]{Theorem} \newtheorem{pro}[df]{Proposition}
\newtheorem{cor}[df]{Corollary} 
\newtheorem{rema}[df] {Remark} 
\def\sfstp{{\hskip-1em}{\bf.}{\hskip1em}}

\def\subject#1{\renewcommand{\thefootnote}{}\footnote
{AMS(MOS) subject classification (2010). Primary: {#1}}}

\def\keywords#1{\renewcommand{\thefootnote}{}\footnote
{Keywords: {#1}}}

\def\enddemo{\qed \endtrivlist} \expandafter\let\csname
enddemo*\endcsname=\enddemo

\def\qedsymbol{\ifmmode\bgroup\else$\bgroup\aftergroup$\fi
\vcenter{\hrule\hbox{\vrule
height.5em\kern.5em\vrule}\hrule}\egroup}
\def\qed{\ifmmode\else\unskip\nobreak\fi\quad\qedsymbol}

\title{\bf Left $m$-invertibility by the adjoint of Drazin inverse and $m$-selfadjointness of Hilbert spaces}
\author{\normalsize B.P.~Duggal, I.H.~Kim}
\date{}

\begin{document}

\maketitle \thispagestyle{empty} \vskip-16pt

\subject{47A05, 47A55; Secondary47A11, 47B47.} \keywords{ Hilbert space,  Left/right multiplication operator,  $m$-left invertible, $m$-isometric and $m$-selfadjoint operators, product of operators, perturbation by nilpotents, commuting operators.  }
\footnote{The work of the second author was supported by the Incheon National University Research Grant in 2017}
\begin{abstract} A Hilbert space operator $A\in\B$ is left $(X,m)$-invertible by $B\in\B$ (resp., $B\in\B$ is an $(X,m)$-adjoint of $A\in\B$) for some operator $X\in\B$ if  $\triangle_{B,A}^m(X)=\sum_{j=0}^m(-1)^j\left(\begin{array}{clcr}m\\j\end{array}\right)B^{m-j}XA^{m-j}=0$  (resp.,
$\delta_{B,A}^m(X)=\sum_{j=0}^m(-1)^j\left(\begin{array}{clcr}m\\j\end{array}\right)B^{(m-j)}XA^j=0$). No Drazin invertible operator $A\in\B$, with Drazin inverse $A_d$, can be left $(I,m)$-invertible (equivalently, $m$-invertible) by its adjoint or its Drazin inverse or the adjoint of its Drazin inverse.
For Drazin inverrtible operators $A$, it is seen that the existence of an $X$ acts as a conduit for implications $\triangle_{B,A}(X)=0\Longrightarrow \delta^m_{C,A}(X)=0$, where  the pair $(B,C)=$ either $(A,A_d)$ or $(A_d,A)$ or $(A^*,A^*_d)$ or $(A^*_d,A^*)$. Reverse implications fail. Assuming certain
commutativity conditions, it is seen that $\triangle_{A^*_d,A}^m(X)=0=\triangle^n_{B^*_d,B}(Y)$ implies $\delta^{m+n-1}_{A^*B^*,AB}(XY)=0=\delta^{m+n-1}_{A^*+B^*,A+B}(XY)$.
\end{abstract}


\section {\sfstp Introduction} Let $\B$ denote the algebra of operators, i.e. bounded linear transformations, on an infinite dimensional complex Hilbert space $\H$ into itself. For $A,B\in\B$, let $L_A$ and $L_B \in B(\B)$ denote respectively the operators
$$
L_A(X)=AX \ {\rm and}\ R_B(X)=XB
$$
of left multiplication by $A$ and right multiplication by $B$. The operator $A$ is left $m$-invertible by $B$, denoted $(B,A)\in$ left-$m$-invertible, if
$$
\triangle_{B,A}^m(I)=\left(L_AR_B-I\right)^m(I)=\sum_{j=0}^m(-1)^j\left(\begin{array}{clcr}m\\j\end{array}\right)B^{m-j}A^{m-j}=0
$$
\cite{DM}. An important class of left $m$-invertible operators, which has been considered by a large number of authors \cite{{AS1},{AS2},{AS3},{OAM},{Fb},{BD1},{BD2},{BD3},{BMN},{BMN1},{BMMN},{G},{G1}}, is that of $m$-isometric operators $A$:
$$
\triangle_{A^*,A}^m(I)=(L_A^*R_A-I)^m(I)=\sum_{j=0}^m(-1)^j\left(\begin{array}{clcr}m\\j\end{array}\right){A^*}^{(m-j)}A^{m-j}=0.
$$
A related, but distinct, class of operators, which has been studied for some time \cite{{Hel},{McR},{TL}}, is that of $m$-selfadjoint operators $A$:
$$
\delta_{A^*,A}^m(I)=(L_A^*-R_A)^m(I)=\sum_{j=0}^m(-1)^j\left(\begin{array}{clcr}m\\j\end{array}\right){A^*}^{(m-j)}A^j=0.
$$

Let $[A,B]=AB-BA$ denote the commutator of $A,B\in \B$. An operator $A\in \B$ is Drazin invertible, with Drazin inverse $A_d$, if
$$
[A_d,A]=0,\ A_d^2A=A_d,\ A^{p+1}A_d=A^p
$$
for some integer $p\geq 1$. (The least integer $p$ for which this holds is then called the Drazin index of $A$.) No Drazin invertible operator $A\in\B$ can be left-$m$-invertible by its adjoint or its Drazin inverse or the adjoint of its Drazin inverse. We say in the following that $A$ is left $(X,m)$-invertible by $B$ (for some operator $X\in\B$), denoted $(B,A)\in$ \ left-$(X,m)$-invertible, if
$$
\triangle_{B,A}^m(X)=\left(L_BR_A-I\right)^m(X)=\sum_{j=0}^m(-1)^j\left(\begin{array}{clcr}m\\j\end{array}\right)B^{m-j}XA^{m-j}=0;
$$
$B\in\B$ is an $(X,m)$-adjoint of $A\in\B$, denoted $(B,A)\in(X,m)$-adjoint, if
$$
\delta_{B,A}^m(X)=(L_B-R_A)^m(X)=\sum_{j=0}^m(-1)^j\left(\begin{array}{clcr}m\\j\end{array}\right)B^{m-j}XA^j=0.
$$
(Here, $A\in (X,m)$-selfadjoint if $(A^*,A)\in (X,m)$-adjoint.) For Drazin invertible operators $A$, there may exist operators $X\in\B$ such that $\triangle_{B,A}^m(X)=0$, where $B=A$ or $A^*$ or $A_d$ or $A_d^*$. We prove that such an operator $X$ necessarily has a representation $X=X_{11}\oplus 0\in B(A^p{(\cal H)}\oplus A^{-p}(0))$ (where $p$ is the Drazin index of $A$). Furthermore, $\triangle_{B,A}^m(X)=0$ implies $\delta_{C,A}^m(X)=0$, where corresponding to $B=A$ or $A^*$ or $A_d$ or $A_d^*$ we have respectively that $C=A_d$ or $A_d^*$ or $A$ or $A^*$. If $A,B\in \B$ are Drazin invertible, (i) $(A_d^*,A)\in$ left-$(X,m)$-invertible and $(B_d^*,B)\in$ left-$(Y,n)$-invertible and (ii) $[A,B]=[X,Y]=[A^*,Y]=[B^*,X]=0$, then $AB$ and $A+B$ are $(XY,m+n-1)$-selfadjoint. The implication $\delta_{A^*+B^*,A+B}^{m+n-1}(XY)=0$ implies $\triangle_{(A^*+B^*)_d,A+B}^{m+n-1}(XY)=0$ fails for Drazin invertible $A,B$ satisfying hypotheses (i) and (ii). Indeed, $A$ and $B$ Drazin invertible does not ensure the Drazin invertibility of $A+B$, even when $A,B$ commute. A sufficient condition in the presence of commutativity is that $AB=0$ \cite{DR}. If the Drazin invertible operators $A,B$ satisfying hypotheses (i) and (ii) satisfy additionally that $AB=0$, then the implication $\delta_{A^*+B^*,A+B}^{m+n-1}(XY)=0$ implies $\triangle_{(A^*+B^*)_d,A+B}^{m+n-1}(XY)=0$ holds. Also, if the Drazin invertible operators $A,B$ satisfy $A\in(X,m)$-isometric, $B\in (Y,n)$-isometric, if (ii) is satisfied and $AB=0$, then  $A+B\in (XY,m+n-1)$-isometric and $\left((A^*+B^*)_d,A+B\right)\in(XY,m+n-1)$-adjoint.

\

\section {\sfstp Results.} Recall, \cite{{FB},{TL}}, that
$$
A\in m-{\rm{isometric}} \Longrightarrow A\in n-{\rm{isometric\ and}}\ A\in m-{\rm selfadjoint}\Longrightarrow  A\in n-{\rm selfadjoint}
$$
for all integers $n\ge m$; again, if $A$ is invertible then
$$
A\in m-{\rm{isometric}} \Longrightarrow A^{-1}\in m-{\rm{isometric\ and}}\ A\in m-{\rm selfadjoint}\Longrightarrow  A^{-1}\in n-{\rm selfadjoint}.
$$

The following proposition says that these results extend to $(X,m)$-operators.

\begin{pro}\label{pro01} Given operators $A,B\in \B$,
\vskip4pt\noindent (i) $(B,A)\in$ left-$(X,m)$-invertible $\Longleftrightarrow (B,A)\in$ left-$(X,n)$-invertible, and
$$
(B,A)\in (X,m)-adjoint \Longleftrightarrow (B,A)\in (X,n)-adjoint
$$
for all integers $n\geq m$;

\vskip4pt\noindent (ii) if $A,B$ are invertible, then
$$
(B,A)\in left-(X,m)-invertible \Longleftrightarrow (B^{-1},A^{-1})\in left-(X,m)-invertible,\ and
$$
$$
(B,A)\in (X,m)-adjoint \Longleftrightarrow (B^{-1},A^{-1})\in (X,m)-adjoint.
$$
\end{pro}

 \begin{demo} (i) The backward implication is evident, and the proof of the forward implication follows from
$$
\triangle_{B,A}^n(X)=\triangle_{B,A}^{n-m}\left(\triangle_{B,A}^m(X)\right)\ {\rm and}\ \delta_{B,A}^n(X)=\delta_{B,A}^{n-m}\left(\delta_{B,A}^m(X)\right)
$$
for all integers $n\geq m$.

\vskip4pt\noindent (ii) If $A,B$ are invertible, then
$$
\triangle_{B^{-1},A^{-1}}^m(X)=L_{B^{-m}}R_{A^{-m}}\left((-1)^m\triangle_{B,A}^m(X)\right)$$
and
$$
 \delta_{B^{-1},A^{-1}}^m(X)=L_{B^{-m}}R_{A^{-m}}\left((-1)^m\delta_{B,A}^m(X)\right).
$$
The proof follows.
\end{demo}

It is well known, \cite[Corollary 2.9]{TL}, that if $A,B\in \B$ are two commutating operators, $A\in m$-selfadjoint and $B\in n$-selfadjoint, then $AB$ and $A+B$ are $(m+n-1)$-selfadjoint; again, if $A\in m$-selfadjoint and $N\in\B$ is a $q$-nilpotent which commutates with $A$, then $A$ is $(m+2q-2)$-selfadjoint. These results extend to  $(X,m)$-selfadjoint operators, as the following proposition proves. The argument we use to prove the proposition differs from most extant proofs (proving similar results); it is similar in spirit to the argument of the proof of Proposition \ref{pro01}, and depends upon a juducious use of some elementary properties of the left and the right multiplication operators.

\begin{pro}\label{pro02} Given operators $X,Y,A,B\in \B$, if :

\vskip4pt\noindent (i) $[A,B]=[X,Y]=[A^*,Y]=[B^*,X]=0$, $A\in (X,m)$-selfadjoint and $B\in (Y,n)$-selfadjoint, then $AB$ and $A+B\in (XY,m+n-1)$-selfadjoint.

\vskip4pt\noindent (ii) $A\in (X,m)$-selfadjoint and $N\in\B$ is a $q$-nilpotent which commutates with $A$, then $\delta_{A^*,A+N}^{m+q-1}(X)=0$. Consequently, $A+N$ is $(m+2q-2)$-selfadjoint.
\end{pro}

\begin{demo} (i) Since the left and right multiplication operators commute, the commutativity hypothesis implies
\begin{eqnarray*}
\delta_{A^*B^*,AB}^{m+n-1}(XY)&=& \left(L_{A^*}L_{B^*}-R_AR_B\right)^{m+n-1}(XY)\\
&=& \left\{(L_{A^*}-R_A)L_{B^*}+R_A(L_{B^*}-R_B)\right\}^{m+n-1}(XY)\\
&=& \left\{\delta_{A^*,A}L_{B^*}+R_A\delta_{B^*,B}\right\}^{m+n-1}(XY)\\
&=& \left\{\sum_{j=0}^{m+n-1}\left(\begin{array}{clcr}m+n-1\\j\end{array}\right)\delta_{A^*,A}^{m+n-1-j}L_{B^*}^{m+n-1-j}R_A^j\delta_{B^*,B}^j\right\}(XY)\\
&=& \sum_{j=0}^{m+n-1}\left(\begin{array}{clcr}m+n-1\\j\end{array}\right)L_{B^*}^{m+n-1-j}R_A^j\left\{\delta_{A^*,A}^{m+n-1-j}\delta_{B^*,B}^j\right\}(XY)\\
&=& \sum_{j=0}^{m+n-1}\left(\begin{array}{clcr}m+n-1\\j\end{array}\right)L_{B^*}^{m+n-1-j}R_A^j\left\{\delta_{B^*,B}^j\delta_{A^*,A}^{m+n-1-j}\right\}(XY)
\end{eqnarray*}
and
\begin{eqnarray*}
\delta_{A^*+B^*,A+B}^{m+n-1}(XY)&=& \left\{\delta_{A^*,A}+\delta_{B^*,B}\right\}^{m+n-1}(XY)\\
&=& \sum_{j=0}^{m+n-1}\left(\begin{array}{clcr}m+n-1\\j\end{array}\right)\left\{\delta_{A^*,A}^{m+n-1-j}\delta_{B^*,B}^j\right\}(XY)\\
&=& \sum_{j=0}^{m+n-1}\left(\begin{array}{clcr}m+n-1\\j\end{array}\right)\left\{\delta_{B^*,B}^j\delta_{A^*,A}^{m+n-1-j}\right\}(XY)
\end{eqnarray*}
The hypothesis  $[X,Y]=[A^*,Y]=[B^*,X]=0$ implies
\begin{eqnarray*}
\delta_{B^*,B}^j(XY)&=&\delta_{B^*,B}^{j-1}(B^*XY-XYB)\\
&=&\delta_{B^*,B}^{j-1}\left\{X(B^*Y-YB)\right\}\\
&=&\delta_{B^*,B}^{j-1}\left\{X\delta_{B^*,B}(Y)\right\}\\
&=&X\delta_{B^*,B}^j(Y)
\end{eqnarray*}
and (similarly)
$$
\delta_{A^*,A}^{m+n-1-j}(XY)=\delta_{A^*,A}^{m+n-2-j}\left\{Y\delta_{A^*,A}(X)\right\}=Y\delta_{A^*,A}^{m+n-1-j}(X).
$$
Hence, since $\delta_{B^*,B}^j(Y)=0$ for all $j\geq n$ and $\delta_{A^*,A}^{m+n-1-j}(X)=0$ for all $m+n-1-j\geq m+n-1-(n-1)=m$, the proof follows.

\vskip4pt\noindent (ii) In this case:
\begin{eqnarray*}
\delta_{A^*,A+N}^{m+q-1}(X)&=& \left\{\delta_{A^*,A}-R_N\right\}^{m+q-1}(X)\\
&=& \sum_{j=0}^{m+q-1}(-1)^j\left(\begin{array}{clcr}m+q-1\\j\end{array}\right)\left\{\delta_{A^*,A}^{m+q-1-j}R_N^j\right\}(X)\\
&=& \sum_{j=0}^{m+q-1}(-1)^j\left(\begin{array}{clcr}m+q-1\\j\end{array}\right)\left\{R_N^j\delta_{A^*,A}^{m+q-1-j}\right\}(X).
\end{eqnarray*}
Since $N^j=0$ for all $j\geq q$ and $\delta_{A^*,A}^{m+q-1-j}(X)=0$ for all $m+q-1-j\geq m+q-1-(q-1)=m$, the proof follows.
\end{demo}

\begin{rema}\label{rema01} Let $\X$ be a complex infinite dimensional Banach space, and let $B(\X)$ denote the algebra of operators on $\X$ into itself. The definition of left-$(X,m)$-invertible operators is equally valid for Banach space operators, and Proposition \ref{pro02} has a Banach space version for products of commuting left-$(X_i,m_i)$-operators, and for perturbation by commuting nilpotents of left-$(X,m)$-invertible operators.
\end{rema}

{\em Given operators $A_i,B_i,X_i\in \B, i=1,2,$ if:
\vskip4pt\noindent  (i) $[A_1,A_2]=[A_1,B_2]=[X_1,X_2]=[A_1,X_2]=[A_2,X_1]=0$ and $(B_i,A_i)\in$ left-$(X_i,m_i)$-invertible, then $(B_1B_2,A_1A_2)\in$ left-$(X_1X_2,m_1+m_2-1)$-invertible;

\vskip4pt\noindent (ii) $N\in \b$ is $q$-nilpotent, $[A_1,N]=0$ and $(B_1,A_1)\in$ left-$(X_1,m_1)$-invertible, then $\triangle_{B_1,A_1+N}^{m_1+q-1}(X_1)=0$ and $(B_1+N,A_1+N)\in$ left-$(m_1+2q-2)$-invertible.}

\vskip6pt\noindent A proof of this follows from the following argument. We have :
\begin{eqnarray*}
& &\triangle_{B_1B_2,A_1A_2}^{m_1+m_2-1}(X_1X_2)\\
&=& \left\{L_{B_2}\triangle_{B_1,A_1}R_{A_2}+\triangle_{B_2,A_2}\right\}^{m_1+m_2-1}(X_1X_2)\\
&=& \sum_{j=0}^{m_1+m_2-1}(-1)^j\left(\begin{array}{clcr}m_1+m_2-1\\j\end{array}\right)\left\{\left(L_{B_2}R_{A_2}\right)^{m_1+m_2-1-j}\triangle_{B_1,A_1}^{m_1+m_2-1-j}\triangle_{B_2,A_2}^j
\right\}(X_1X_2),
\end{eqnarray*}
where $\left[\triangle_{B_1,A_1},\triangle_{B_2,A_2}\right]=0$, $\triangle_{B_2,A_2}^j(X_2)=0$ for all $j\geq m_2$ and $\triangle_{B_1,A_1}^{m_1+m_2-1-j}(X_1)=0$ for all $j\leq m_2-1$;
\begin{eqnarray*}
\triangle_{B_1,A_1+N}^{m_1+q-1}(X_1)&=& \left\{\triangle_{B_1,A_1}+L_{B_1}R_N\right\}^{m_1+q-1}(X_1)\\
&=& \sum_{j=0}^{m_1+q-1}\left(\begin{array}{clcr}m_1+q-1\\j\end{array}\right)\left\{\triangle_{B_1,A_1}^{m_1+q-1-j}\left(L_{B_1}R_N\right)^j\right\}(X_1),
\end{eqnarray*}
where $\left[\triangle_{B_1,A_1},L_{B_1}R_N\right]=0$, $R_N^j=0$ for all $j\geq q$ and $\triangle_{B_1,A_1}^{m_1+q-1-j}(X_1)=0$ for all $j\leq q-1$

The following corollary is immediate from the argument of the proof of Proposition \ref{pro02}.

\begin{cor}\label{cor01} If $A,B\in\B$ are commuting operators such that $A\in (X,m)$-selfadjoint and $B\in (X,n)$-selfadjoint for some operator $X\in\B$, then $AB$ and $A+B$ are $(X,m+n-1)$-selfadjoint.
\end{cor}

\vskip4pt\noindent {\bf Drazin invertible operators $A$.} The ascent of an operator $A$, asc$(A)$, (resp., the descent of an operator $A$, dsc$(A)$) is the least positive integer $p$ such that $A^{-p}(0)=A^{-(p+1)}(0)$ (resp., $A^p(\H)=A^{p+1}(\H)$). The Drazin invertibility ensures the existence of such an integers $p$, and then ${\rm asc}(A)={\rm dsc}(A)=p$; the integer $p$ is the Drazin index of $A$ \cite{DR}: Throughout the following, unless otherwise stated, we assume that $A$ is Drazin invertible and that the Drazin index of $A$ is $p$.

The Drazin invertibility of $A$ induces a decomposition
$$
\H=A^p(\H)\oplus A^{-p}(0)=\H_1 \oplus \H_2
$$
of $\H$, and a decomposition
$$
A=A\mid_{\H_1}\oplus A\mid_{\H_2}=A_1\oplus A_2, \ A_1 \ {\rm invertible \ and }\ A_2 \ p{\rm -nilpotent},
$$
of $A$ \cite{DR}. Accordingly, the Drazin inverse $A_d$ of $A$ has a decomposition
$$
A_d=A_1^{-1}\oplus 0\in B(\H_1\oplus\H_2)
$$
\cite[Theorem 2.23]{DR}. A Drazin invertible operator can not be $m$-isometric (reason: if it is then the spectrum $\sigma(A)$ of $A$ is the union of a subset of the boundary $\partial\mathbb D$ of the unit disc $\mathbb D$ with the point set $\{0\}$ and the spectrum of an $m$-isometry is either the closure $\overline{\mathbb D}$ of the unit disc $\mathbb D$ or a subset of $\partial\mathbb D$ \cite{FB}). Again, a Drazin invertible operator $A$ can not be left-$m$-invertible by $A_d$ or $A_d^*$: this is consequent from the fact that
$$
\sum_{j=0}^m(-1)^m\left(\begin{array}{clcr}m\\j\end{array}\right)B^{m-j}A^{m-j}=\left\{\sum_{j=0}^m(-1)^m\left(\begin{array}{clcr}m\\j\end{array}\right)B_1^{-m+j}A_1^{m-j}
\right\}\oplus I_2 \neq 0,
$$
where $I_2=I\mid_{\H_2}$, $B=A_d$ or $A_d^*$ and $B_1=A_1$ or (resp.,) $A_1^*$. There may, however, exist operators $X\in\B$ such that $A$ is left-$(X,m)$-invertible by $A^*$ or $A_d$ or $A_d^*$ or (even) $A$. For example, if $\sum_{j=0}^m(-1)^m\left(\begin{array}{clcr}m\\j\end{array}\right)A_1^{-m+j}X_{11}A_1^{m-j}=0$ (resp., $\sum_{j=0}^m(-1)^m\left(\begin{array}{clcr}m\\j\end{array}\right){A_1^*}^{-m+j}X_{11}A_1^{m-j}=0$) for some $X_{11}\in B(\H_1)$, then $(A_d,A)\in$ left-$(X,m)$-invertible (resp., $(A_d^*,A)\in$ left-$(X,m)$-invertible) for $X=X_{11}\oplus 0\in B(\H_1\oplus \H_2)$.

\begin{rema}\label{rema02} In contrast to the situation for left $m$-invertible operators, every operator $A$ satisfies $\delta_{A,A}^m(I)=0$ for all integers $m\geq 1$. If $A$ is Drazin invertible and $m$-selfadjoint, then
$$
\delta_{A^*,A}^m(I)=0 \Longleftrightarrow \sum_{j=0}^m(-1)^j\left(\begin{array}{clcr}m\\j\end{array}\right)\left\{{A_1^*}^{m-j}A_1^j\oplus {A_2^*}^{m-j}A_2^j\right\}=0,
$$
i.e., $A$ is $m$-selfadjoint if and only if the invertible operator $A_1$ and the $p$-nilpotent operator $A_2$ are $m$-selfadjoint. In particular, if $m=2$, then $A$ is 2-selfadjoint if and only if $A=A_1\oplus 0$ (see \cite[Theorem 3.1]{McR}, where it is proved that a $\B$ operator is 2-selfadjoint if and only if it is selfadjoint). If $\delta_{A_d^*,A}^m(I)=0$, then
\begin{eqnarray*}
\sum_{j=0}^m(-1)^j\left(\begin{array}{clcr}m\\j\end{array}\right){A_1^*}^{-m+j}A_1^j=0 & \Longleftrightarrow& \sum_{j=0}^m(-1)^j\left(\begin{array}{clcr}m\\j\end{array}\right){A_1^*}^jA_1^j=0\\ &\Longleftrightarrow& \triangle^m_{{A_1^*}^{-1}A_1}(I)=0
\end{eqnarray*}
(where $I_1=I\mid_{\H_1}$). Consequently, $\sigma(A)$  is the union of a subset of $\partial\mathbb D$ with the point set $\{0\}$. In the particular case in which $m=2$, this implies (either $A=0$, or) $\sigma(A)=\{\pm{1},0\}$.
\end{rema}

The following theorem proves that if $(B,A)\in$ left-$(X,m)$-invertible for some $X\in\B$, and $B=A$ or $A^*$ or $A_d$ or $A_d^*$, then necessarily $X=X_{11}\oplus 0\in B(\H_1\oplus\H_2)$. As a consequence it is seen that if $A$ is left-$(X,m)$-invertible by $A_d^*$, then $A$ is $(X,m)$-selfadjoint.

\begin{thm}\label{thm01} Let $A,X\in\B$, where $A$ is Drazin invertible with Drazin inverse $A_d$.

\noindent (i) If $(B,A)\in$ left-$(X,m)$-invertible, and $B=A$ or $A^*$ or $A_d$ or $A_d^*$, then $X=X_{11}\oplus 0\in B(\H_1\oplus\H_2)$.

Consequently, one has:

\vskip4pt\noindent (ii)
\begin{eqnarray*}
\triangle_{A,A}^m(X)=0 &\Longrightarrow& \delta_{A_d,A}^m(X)=0;\\
\triangle_{A_d,A}^m(X)=0 &\Longrightarrow& \delta_{A,A}^m(X)=0;\\
\triangle_{A^*,A}^m(X)=0 &\Longrightarrow& \delta_{A_d^*,A}^m(X)=0;\\
\triangle_{A_d^*,A}^m(X)=0 &\Longrightarrow& \delta_{A^*,A}^m(X)=0.
\end{eqnarray*}
 \end{thm}

\begin{demo} The proof for all four choices of $B$ is similar: we consider the case in which $B=A$. Letting, as above, $A=A_1\oplus A_2\in B(\H_1\oplus\H_2)$, where $A_1$ is invertible and $A_2$ is $p$-nilpotent, and letting $X\in B(\H_1\oplus\H_2)$ have the matrix representation $X=\left[X_{ik}\right]_{i,k=1}^2$, we have:
\begin{eqnarray*}
\triangle^m_{A,A}(X) &=& \sum_{j=0}^m(-1)^j\left(\begin{array}{clcr}m\\j\end{array}\right){A}^{m-j}XA^{m-j}=0 \\ &\Longleftrightarrow& \left[\sum_{j=0}^m(-1)^j\left(\begin{array}{clcr}m\\j\end{array}\right){A_i}^{m-j}X_{ik}A_k^{m-j}\right]_{i,k=1}^2=0
\end{eqnarray*}
(Here, the $(2,1)$ and $(2,2)$ entries equal $X_{21}$ and $X_{22}$, respectively, in the case in which $B=A_d$ or $A_d^*$). We have two possibilities: either $p<m$ or $p\geq m$. If $p\geq m$, then Proposition \ref{pro01} tell us that $(A,A)\in$ left-$(X,n)$-invertible for all integers $n\geq p$. Hence it will suffice to prove $X_{12}=X_{21}=X_{22}=0$ for $p<m$. We consider the case of $X_{12}$: the proof for the other two cases is similar. If $p<m$, then $A_1^tX_{12}A_2^t=0$ for all $t\geq p$. We have:
\begin{eqnarray*}
& & \sum_{j=0}^m(-1)^j\left(\begin{array}{clcr}m\\j\end{array}\right){A}_1^{m-j}X_{12}A_2^{m-j}=0\\  &\Longrightarrow& A_1^{p-1}\left\{\sum_{j=0}^m(-1)^j\left(\begin{array}{clcr}m\\j\end{array}\right){A_1}^{m-j}X_{12}A_2^{m-j}\right\}A_2^{p-1}=0\\
&\Longleftrightarrow& {A_1}^{p-1}X_{12}A_2^{p-1}=0\\
&\Longrightarrow& A_1^tX_{12}A_2^t=0 \ {\rm for \ all}\ t\geq p- 1
\end{eqnarray*}
Repeating this argument a further $(p-2)$-times we obtain $A_1^tX_{12}A_2^t=0$ for all $t\geq 1$. Hence
$$
0=\sum_{j=0}^m(-1)^j\left(\begin{array}{clcr}m\\j\end{array}\right){A}_1^{m-j}X_{12}A_2^{m-j}=X_{12}.
$$

\noindent (ii) The proof of (i) implies
\begin{eqnarray*}
& & \triangle_{A,A}^m(X)=0\\ &\Longleftrightarrow& \sum_{j=0}^m(-1)^j\left(\begin{array}{clcr}m\\j\end{array}\right){A}_1^{m-j}X_{11}A_1^{m-j}=0,\ X_{12}=X_{21}=X_{22}=0\\
&\Longleftrightarrow& \sum_{j=0}^m(-1)^j\left(\begin{array}{clcr}m\\j\end{array}\right){A}_1^{-m+j}X_{11}A_1^{-m+j}=0,\ X_{12}=X_{21}=X_{22}=0\\
& & \ \ \ \ \ \ \ \ ({\rm since} \ A_1 \ {\rm is \ invertible})\\
&\Longleftrightarrow& \sum_{j=0}^m(-1)^j\left(\begin{array}{clcr}m\\j\end{array}\right){A}_1^{-(m-j)}X_{11}A_1^j=0,\ X_{12}=X_{21}=X_{22}=0\\
&\Longrightarrow& \delta_{A_d,A}^m(X)=0;
\end{eqnarray*}
\begin{eqnarray*}
& & \triangle_{A_d,A}^m(X)=0\\  &\Longleftrightarrow& \sum_{j=0}^m(-1)^j\left(\begin{array}{clcr}m\\j\end{array}\right){A}_1^{-m+j}X_{11}A_1^{m-j}=0,\ X_{12}=X_{21}=X_{22}=0\\
&\Longleftrightarrow& \sum_{j=0}^m(-1)^j\left(\begin{array}{clcr}m\\j\end{array}\right){A}_1^{m-j}X_{11}A_1^{-m+j}=0,\ X_{12}=X_{21}=X_{22}=0\\
& & \ \ \ \ \ \ \ \ ({\rm since} \ A_1 \ {\rm is \ invertible})\\
&\Longleftrightarrow& \sum_{j=0}^m(-1)^j\left(\begin{array}{clcr}m\\j\end{array}\right){A}_1^{m-j}X_{11}A_1^j=0,\ X_{12}=X_{21}=X_{22}=0\\
&\Longrightarrow& \delta_{A,A}^m(X)=0;
\end{eqnarray*}
\begin{eqnarray*}
& & \triangle_{A^*,A}^m(X)=0\\ &\Longleftrightarrow& \sum_{j=0}^m(-1)^j\left(\begin{array}{clcr}m\\j\end{array}\right){A^*}_1^{(m-j)}X_{11}A_1^{m-j}=0,\ X_{12}=X_{21}=X_{22}=0\\
&\Longleftrightarrow& \sum_{j=0}^m(-1)^j\left(\begin{array}{clcr}m\\j\end{array}\right){A^*}_1^{(-m+j)}X_{11}A_1^{-m+j}=0,\ X_{12}=X_{21}=X_{22}=0\\
& & \ \ \ \ \ \ \ \ ({\rm since} \ A_1 \ {\rm is \ invertible})\\
&\Longleftrightarrow& \sum_{j=0}^m(-1)^j\left(\begin{array}{clcr}m\\j\end{array}\right){A^*}_1^{(-m+j)}X_{11}A_1^j=0,\ X_{12}=X_{21}=X_{22}=0\\
&\Longrightarrow& \delta_{A_d^*,A}^m(X)=0
\end{eqnarray*}
and finally
\begin{eqnarray*}
& &\triangle_{A_d^*,A}^m(X)=0\\ &\Longleftrightarrow& \sum_{j=0}^m(-1)^j\left(\begin{array}{clcr}m\\j\end{array}\right){A^*}_1^{(-m+j)}X_{11}A_1^{m-j}=0,\ X_{12}=X_{21}=X_{22}=0\\
&\Longleftrightarrow& \sum_{j=0}^m(-1)^j\left(\begin{array}{clcr}m\\j\end{array}\right){A^*}_1^{m-j}X_{11}A_1^{-m+j}=0,\ X_{12}=X_{21}=X_{22}=0\\
& & \ \ \ \ \ \ \ \ ({\rm since} \ A_1 \ {\rm is \ invertible})\\
&\Longleftrightarrow& \sum_{j=0}^m(-1)^j\left(\begin{array}{clcr}m\\j\end{array}\right){A^*}_1^{m-j}X_{11}A_1^j=0,\ X_{12}=X_{21}=X_{22}=0\\
&\Longrightarrow& \delta_{A^*,A}^m(X)=0.
\end{eqnarray*}
\end{demo}

\begin{rema}\label{rema03} The reverse implications in Theorem \ref{thm01}(ii) fail. This is for the reason that $\delta_{B,A}^m(X)=0$, $B=A$ or $A^*$ or $A_d$ or $A_d^*$, does not imply $X=X_{11}\oplus 0$. Consider, for example, $\delta_{A,A}^m(X)=0$. Assuming, without loss of generality, that $p<m$, it is seen that
$$
\delta_{A,A}^m(X)=0 \Longleftrightarrow \left[\sum_{j=0}^m(-1)^j\left(\begin{array}{clcr}m\\j\end{array}\right){A_i}^{m-j}X_{ik}A_k^j\right]_{i,k=1}^2=0;
$$
\begin{eqnarray*}
& & \sum_{j=0}^m(-1)^j\left(\begin{array}{clcr}m\\j\end{array}\right){A}_1^{m-j}X_{12}A_2^j=0\\
&\Longrightarrow& \left\{\sum_{j=0}^m(-1)^j\left(\begin{array}{clcr}m\\j\end{array}\right){A}_1^{m-j}X_{12}A_2^j\right\}A_2^{p-1}=0\\
&\Longrightarrow& A_1^mX_{12}A_2^{p-1}=0 \Longleftrightarrow X_{12}A_2^{p-1}=0\\
&\Longrightarrow& \left\{\sum_{j=0}^m(-1)^j\left(\begin{array}{clcr}m\\j\end{array}\right){A}_1^{m-j}X_{12}A_2^j\right\}A_2^{p-2}=0\\
&\Longrightarrow& A_1^mX_{12}A_2^{p-2}=0 \Longleftrightarrow X_{12}A_2^{p-2}=0\\
& & \ \ \ \ \ \ \ \ \cdots\\
&\Longrightarrow& \left\{\sum_{j=0}^m(-1)^j\left(\begin{array}{clcr}m\\j\end{array}\right){A}_1^{m-j}X_{12}A_2^j\right\}A_2=0\\
&\Longrightarrow& A_1^mX_{12}A_2=0\Longleftrightarrow X_{12}A_2=0\\
&\Longrightarrow& A_1^mX_{12}=0 \Longleftrightarrow X_{12}=0.
\end{eqnarray*}
\end{rema}
Similarly, $X_{21}=0$, and hence $\delta_{A,A}^m(X)=0 \Longrightarrow X= X_{11}\oplus X_{22}$. Similar arguments show that $X_{12}=X_{21}=0$ if $B=A^*$, and $X_{12}=0$ if $B=A_d$ or $A_d^*$. Examples proving that the reverse implications fail are not difficult to construct. For example, if $A_1\in B(\H_1)$ is an invertible operator such that $A_1^*X_{11}=X_{11}A_1$ for some $X_{11}\in B(\H_1)$, and $N_2\in B(\H_2)$ is a $2$-nilpotent operator, then $\delta_{A^*,A}^3(X)=0$ for $A=A_1\oplus N_1$ and $X=X_{11}\oplus I_2$. However, $\triangle_{A_d^*,A}^3(X)\neq 0$.

Theorem \ref{thm01} and Remark \ref{rema03} taken together imply that {\em a Drazin invertible operator $A\in\B$ satisfies  $$(A_d^*,A)\in \ {\rm left}-(X,m)-{\rm invertible} \ \Longleftrightarrow A\in (X,m)-{\rm selfadjoint}$$ if and only if $X\in B(\H_1\oplus \H_2)$ has a representation $X=X_{11}\oplus 0$.} The following theorem gives a sufficient condition for $(A_d^*,A)\in$ left-$(X,m)$-invertible to imply $A$ is $n$-selfadjoint (for $n=m+2p-2$).

\begin{thm}\label{thm02} If $\triangle_{A_d^*,A}^m(X)=0$, $[A,X]=0$ and $X_{11}=X\mid_{\H_1}$ has a dense range, then $A$ is $(m+2p-2)$-selfadjoint. Furthermore, if $m=2$, then $A$ is a $(2p-1)$-selfadjoint operator (which is a perturbation by a $p$-nilpotent operator of a selfadjoint operator).
\end{thm}

\begin{demo} If $\triangle_{A_d^*,A}^m(X)=0$, then $X=X_{11}\oplus 0\in B(\H_1\oplus \H_2)$, and this if $[X,A]=0$ implies $[X_{11},A_1]=0$. Hence
\begin{eqnarray*}
\triangle_{A_d^*,A}^m(X)=0 &\Longleftrightarrow& \sum_{j=0}^m(-1)^j\left(\begin{array}{clcr}m\\j\end{array}\right){A^*}_1^{(-m+j)}X_{11}A_1^{m-j}=0\\
&\Longleftrightarrow& \sum_{j=0}^m(-1)^j\left(\begin{array}{clcr}m\\j\end{array}\right){A^*}_1^{(m-j)}X_{11}A_1^{-m+j}=0\\
& & (\rm{since} \ A_1  \ \rm{is \ invertible})\\
&\Longleftrightarrow& \sum_{j=0}^m(-1)^j\left(\begin{array}{clcr}m\\j\end{array}\right){A^*}_1^{(m-j)}X_{11}A_1^{j}=0\\
&\Longleftrightarrow& \left\{\sum_{j=0}^m(-1)^j\left(\begin{array}{clcr}m\\j\end{array}\right){A^*}_1^{(m-j)}A_1^j\right\}X_{11}=0.
\end{eqnarray*}
This, if $X_{11}$ has a dense range, implies
$$
\sum_{j=0}^m(-1)^j\left(\begin{array}{clcr}m\\j\end{array}\right){A^*}_1^{(m-j)}A_1^j=0,
$$
i.e., $A_1\in B(\H_1)$ is $m$-selfadjoint. Consequently, $A_1\oplus 0\in B(\H_1\oplus \H_2)$ is $m$-selfadjoint. Now define the $p$-nilpotent operator $N$ by  $N=0\oplus A_2\in B(\H_1\oplus \H_2)$. Then $N$ commutes with $A_1\oplus 0$, and hence $A=(A_1\oplus 0)+N$ is $(m+2p-2)$-selfadjoint.

Now let $m=2$. Then the $2$-selfadjoint operator $A_1$ above is selfadjoint \cite[Theorem 3.1]{McR}, and the operator $A$ being the perturbation of a selfadjoint operator by a $p$-nilpotent operator is a $(2p-1)$-selfadjoint operator.

\end{demo}

If $A,B\in\B$ are Drazin invertible (with Drazin inverses $A_d, B_d$), $(A_d^*,A)\in$ left-$(X,m)$-invertible, $(B_d^*,A)\in$ left-$(Y,n)$-invertible (for some operators $X,Y\in\B$) and
\begin{eqnarray}
[X,Y]=[A,B]=[A^*,Y]=[B^*,X]=0,
\end{eqnarray}
then
\begin{eqnarray*}
(AB)_d=A_dB_d=B_dA_d, \ [A_d^*,Y]=[B_d^*,Y]=0,\ \triangle_{A_d^*B_d^*,AB}^{m+n-1}(XY)=0
\end{eqnarray*}
(see Remark \ref{rema01}), and hence
 $$
 \delta_{A^*B^*,AB}^{m+n-1}(XY)=0
 $$
 (see Theorem \ref{thm01}). The hypotheses on $A,B$ (in the above) are not sufficient to guarantee the Drazin invertibility of $A+B$. Even given $A,B$ and $A+B$ are all Drazin invertible, do the above hypotheses on $A,B,X$ and $Y$ imply $\triangle_{(A^*+B^*)_d,A+B}^{m+n-1}(XY)=0$ ? Postponing an answer to this question for the time being, we prove in the following theorem that the hypotheses above are sufficient to guarantee $\triangle_{(A^*+B^*)_d,A+B}^{m+n-1}(XY)=0$.

\begin{thm}\label{thm03} If $(A_d^*,A)\in$ left-$(X,m)$-invertible, $(B_d^*,B)\in$ left-$(Y,n)$-invertible and hypothesis (1) holds, then $AB$ and $A+B$ are $(XY,m+n-1)$-selfadjoint.
\end{thm}

\begin{demo}
The hypotheses imply that the implications
$$
\triangle_{A_d^*,A}^m(X)=0 \Longrightarrow \delta_{A^*,A}^m(X)=0 \ {\rm and}\ \triangle_{B_d^*,B}^n(Y)=0 \Longrightarrow \delta_{B^*,B}^n(X)=0
$$
hold (see Theorem \ref{thm01}). Since $A,B,X$ and $Y$ satisfy the hypotheses of Proposition \ref{pro02}, we conclude that $AB$ and $A+B$ are $(XY,m+n-1)$-selfadjoint.
\end{demo}

The Drazin invertibility of $A$ and $B$ does not ensure the Drazin invertibility of $A+B$. (Hence the implication $\delta_{A^*+B^*,A+B}^{m+n-1}(XY)=0 \Longrightarrow \triangle_{(A^*+B^*)_d,A+B}^{m+n-1}(XY)=0$ is not guaranteed by the hypotheses of Theorem \ref{thm03}). A sufficient condition for $A+B$ to be Drazin invertible for commuting Drazin invertible $A$ and $B$ is that $AB=0$ \cite[page 93]{DR}. Assume, as above, that $A=A_1\oplus A_2\in B(\H_1\oplus \H_2)$, $[A,B]=0=AB$ and let $B\in B(\H_1\oplus \H_2)$ have the matrix representation $B=\left[B_{ij}\right]_{i,j=1}^2.$ Then
$$
[A,B]=0 \Longrightarrow B=B_{11}\oplus B_{22}$$
 and if also  $ AB=0$, then $$ B_{11}=0 \ {\rm and \ hence}\ B=0\oplus B_{22},$$

\noindent  $[A_2,B_{22}]=0=A_2B_{22}$, the operators $B_{22}$ and $A_2+B_{22}$ are Drazin invertible and
$$
(A+B)_d=\left(\begin{array}{cc}A_1^{-1}&0\\0 & (A_2+B_{22})_d\end{array}\right)\in B(\H_1\oplus\H_2).
$$

The following theorem answers the question posed above for the case in which $AB=0$.

\begin{thm}\label{thm04} If in addition to the hypotheses of Theorem \ref{thm03}, $A$ and $B$ satisfy the hypothesis that $AB=0$, then
$$
\delta_{(A+B)_d^*,A+B}^{m+n-1}(XY)=0.
$$
\end{thm}

\begin{demo}
If we let $A=A_1\oplus A_2\in B(\H_1\oplus \H_2)$, $A_1$ invertible and $A_2$ $p$-nilpotent, then upon letting $X=\left[X_{ij}\right]_{i,j=1}^2\in B(\H_1\oplus \H_2)$ the hypothesis $(A_d^*,A)\in$ left-$(X,m)$-invertible implies $X=X_{11}\oplus 0$ (see Theorem \ref{thm01}(i)). Let $Y\in B(\H_1\oplus \H_2)$ have the matrix representation $Y=\left[Y_{ij}\right]_{i,j=1}^2\in B(\H_1\oplus \H_2)$. Then $[X,Y]=0$ implies
$$
X_{11}Y_{12}=0=Y_{21}X_{11},
$$
and hence
$$
XY=X_{11}Y_{11}\oplus 0\in B(\H_1\oplus \H_2).
$$
As seen above, $A+B=A_1\oplus (A_2+B_{22})\in B(\H_1\oplus \H_2)$. Hence (see Theorem \ref{thm03})
$$
\delta_{A^*+B^*,A+B}^{m+n-1}(XY)=0 \Longleftrightarrow \left(\delta_{A_1^*,A_1}^{m+n-1}(X_{11}Y_{11})\right)\oplus 0=0.
$$
Since $A_1$ is invertible,
\begin{eqnarray*}
& &\delta_{A_1^*,A_1}^{m+n-1}(X_{11}Y_{11})=0\\
 &\Longleftrightarrow& \delta_{{A_1^*}^{-1},A_1^{-1}}^{m+n-1}(X_{11}Y_{11})=0\\
&\Longleftrightarrow& \sum_{j=0}^{m+n-1}(-1)^j\left(\begin{array}{clcr}m+n-1\\j\end{array}\right){A_1^*}^{(-m-n+1+j)}X_{11}Y_{11}{A_1}^{-j}=0\\
&\Longleftrightarrow& \sum_{j=0}^{m+n-1}(-1)^j\left(\begin{array}{clcr}m+n-1\\j\end{array}\right){A_1^*}^{(-m-n+1+j)}X_{11}Y_{11}{A_1}^{m+n-1-j}=0\\
&\Longleftrightarrow& \triangle_{{A_1^*}^{-1},A_1}^{m+n-1}(X_{11}Y_{11})=0.
\end{eqnarray*}
Since
$$
\triangle_{(A^*+B^*)_d,A+B}^{m+n-1}(XY)=\triangle_{{A_1^*}^{-1},A_1}^{m+n-1}(X_{11}Y_{11})\oplus 0,
$$
we have $$\delta_{(A^*+B^*)_d,A+B}^{m+n-1}(XY)=0.$$
\end{demo}

The following theorem is an analogue of Theorem \ref{thm03} for operators $A\in (X,m)$-isometric, $\triangle_{A^*,A}^m(X)= \sum_{j=0}^m(-1)^j\left(\begin{array}{clcr}m\\j\end{array}\right){A^*}^{(m-j)}XA^{m-j}$. Recall that if $A\in (X,m)$-isometric, then $A\in (X,n)$-isometric for all $n\geq m$.

\begin{thm}\label{thm05} If $A\in (X,m)$-isometric, $B\in (Y,n)$-isometric and
$$
[X,Y]=[A,B]=[A^*,Y]=[B^*,X]=0=AB,
$$
then $A+B$ is $(m+n-1,XY)$-isometric and $\left((A+B)_d^*,A+B\right)\in(m+n-1,XY)$-adjoint.
\end{thm}

\begin{demo}
If $A\in (X,m)$-isometric, then $X=X_{11}\oplus 0\in B(\H_1\oplus \H_2)$ (see Theorem \ref{thm01}) and (as seen above)
$$
A+B=A_1\oplus (A_2+B_{22}),\ (A+B)_d=A_1^{-1}\oplus(A_2+B_{22})_d,\ XY=X_{11}Y_{11}\oplus 0,\ [A_1^*,Y_{11}]=0.
$$
Since $A\in (X,m)$-isometric if and only if $\triangle_{A^*,A}^m(X)=0$, we have
\begin{eqnarray*}
\triangle_{A^*,A}^m(X)=0 &\Longleftrightarrow& \triangle_{A_1^*,A_1}^m(X_{11})=0\\
&\Longleftrightarrow& \triangle_{{A_1^*}^{-1},A_1^{-1}}^m(X_{11})=0 \ \ ({\rm since} \ A_1 \ {\rm is \ invertible})\\
&\Longrightarrow& \triangle_{{A_1^*}^{-1},A_1^{-1}}^{m+n-1}(X_{11})=0  \ {\rm for \ all} \ n\geq 1\\ &\Longleftrightarrow& \left(\triangle_{{A_1^*}^{-1},A_1^{-1}}^{m+n-1}(X_{11})\right)A_1^{m+n-1}=0\\
&\Longrightarrow& \delta_{{A_1^*}^{-1},A_1}^{m+n-1}(X_{11})=0\\
& \Longrightarrow& \delta_{{A_1^*}^{-1},A_1}^{m+n-1}(X_{11}Y_{11})=0,
\end{eqnarray*}
where the final implication follows from the fact that
\begin{eqnarray*}
\delta_{{A_1^*}^{-1},A_1}(X_{11}Y_{11})&=&{A_1^*}^{-1}X_{11}Y_{11}-X_{11}Y_{11}A_1\\
&=&\left({A_1^*}^{-1}Y_{11}\right)X_{11}-Y_{11}X_{11}A_1\\
&=&Y_{11}\left({A_1^*}^{-1}X_{11}-X_{11}A_1\right)\\
&=&Y_{11}\delta_{{A_1^*}^{-1},A_1}(X_{11}).
\end{eqnarray*}
Since
$$
\delta_{{A_1^*}^{-1},A_1}^{m+n-1}(X_{11}Y_{11})=0 \Longleftrightarrow \delta_{(A+B)_d^*,A+B}^{m+n-1}(XY)=0,
$$
we conclude $\left((A+B)_d^*,A+B\right)\in (XY, m+n-1)$-adjoint. To complete the proof, we observe that
\begin{eqnarray*}
A\in (X,m)-{\rm isometric} &\Longrightarrow&  \triangle_{A_1^*,A_1}^{m+n-1}(X_{11})=0\\
&\Longrightarrow& \triangle_{A_1^*,A_1}^{m+n-1}(X_{11}Y_{11})=0\\
&\Longleftrightarrow& \triangle_{(A+B)^*,A+B}^{m+n-1}(XY)=0.
\end{eqnarray*}
\end{demo}


\vskip10pt \noindent\normalsize\rm B.P. Duggal,{University of Ni\v s,
Faculty of Sciences and Mathematics,
P.O. Box 224, 18000 Ni\v s, Serbia}.

\noindent\normalsize \tt e-mail:  bpduggal@yahoo.co.uk

\vskip6pt\noindent \noindent\normalsize\rm I. H. Kim, Department of
Mathematics, Incheon National University, Incheon, 22012, Korea.\\
\noindent\normalsize \tt e-mail: ihkim@inu.ac.kr

\end{document}